\documentclass[10pt]{article}
\usepackage{amsmath,amssymb,latexsym}
\usepackage{theorem}
\usepackage{color}
\usepackage{comment}
\newtheorem{theorem}{Theorem}
\newtheorem{lemma}{Lemma}

\newtheorem{corollary}{Corollary}
\newtheorem{proposition}{Proposition}
{\theorembodyfont{\rmfamily} 
\newtheorem{remark}{Remark}}
{\theorembodyfont{\slshape} }

\newcommand{\field}[1]{\mathbb{#1}}
\newcommand{\R}{\field{R}}

\newcommand{\N}{\field{N}}
\newcommand{\C}{\field{C}}

\newcommand{\BB}{{\mathcal B}}
\renewcommand{\SS}{{\mathcal S}}

\newcommand{\EE}{{\mathcal E}}
\newcommand{\TT}{{\mathcal T}}
\newcommand{\UU}{{\mathcal U}}

\renewcommand{\Re}{\mathop{\rm Re}}

\DeclareRobustCommand{\qed}{%
\ifmmode 
\else \leavevmode\unskip\penalty9999 \hbox{}\nobreak\hfill \fi
\quad\hbox{\qedsymbol}}
\newcommand{\openbox}{\leavevmode
\hbox to.77778em{%
\hfil\vrule
\vbox to.675em{\hrule width.6em\vfil\hrule}%
\vrule\hfil}}
\newcommand{\qedsymbol}{\openbox}
\newcommand{\proofname}{Proof}
\newenvironment{proof}[1][\proofname]{\par
\normalfont \trivlist \item[\hskip\labelsep   \itshape #1. ]
\ignorespaces
}{%
\qed\endtrivlist }
\title{Asymptotic upper bounds for the entropy of orthogonal
polynomials in the Szeg\H{o} class}

\author{ B.\ Beckermann\\ Department of Mathematics, Universite  de Lille I,
France  \and A.\ Mart{\'\i}nez-Finkelshtein\thanks{Corresponding author.
E-mail: \texttt{andrei@ual.es}}\\
Department of Statistics and Applied Mathematics, University of
Almer{\'\i}a \\ and Institute Carlos I for Theoretical and
Computational
Physics, \\ Granada University, Spain \and E.A.\ Rakhmanov\\
Department of Mathematics, University of South Florida, USA \and
F.\ Wielonsky\\ Department of Mathematics, Universite  de Lille I,
 \\ and INRIA Sophia Antipolis, France}

\date{}

\begin{document}

\maketitle
\begin{abstract}
We give an asymptotic upper bound as $n\to\infty$ for the entropy
integral
$$E_n(w)= -\int p_n^2(x)\log (p_n^2(x))w(x)dx,$$ where $p_n$ is the $n$th
degree orthonormal polynomial with respect to a weight $w(x)$ on
$[-1,1]$ which belongs to the Szeg\H{o} class. We also study two
functionals closely related to the entropy integral. First, their
asymptotic behavior is completely described for weights $w$ in the
Bernstein class. Then, as for the entropy, we obtain asymptotic
upper bounds for these two functionals when $w(x)$ belongs to the
Szeg\H{o} class. In each case, we give conditions for these upper
bounds to be attained.
\end{abstract}

\section{Introduction}

In the framework of the density functional theory (see e.g.\
\cite{dg90,py89}) the physical and chemical properties of
fermionic systems are described by means of the single-particle
probability densities. If $\Psi(\vec{r})$ is the wave function of
a single-particle system in a ($D$-dimensional) position space,
and $\widehat{\Psi}(\vec{p})$ is the corresponding wave function
in momentum space (that is, the Fourier transform of
$\Psi(\vec{r})$), then the  position and momentum densities of the
system are given by
\[
\rho(\vec{r}) = | \Psi(\vec{r}) |^2 \;, \qquad \gamma(\vec{p}) = |
\widehat{\Psi}(\vec{p}) |^2 \,,
\]
respectively. It is known that the Boltzmann-Gibbs-Shannon
position-space entropy
\begin{equation*} \label{eserho}
S(\rho) = -\int \rho(\vec{r}) \ln \rho(\vec{r}) \, d \vec{r}
\end{equation*}
measures the uncertainty in the localization of  the particle in
space (lower entropy indicates a more concentrated wave function,
with the associated higher accuracy in predicting the localization
of the particle). Similar is true for the momentum-space entropy
\begin{equation*} \label{esegamma}
S(\gamma) = -\int \gamma(\vec{p}) \ln \gamma(\vec{p}) \, d
\vec{p}\,.
\end{equation*}
These quantities have importance in the study of the structure and
dynamics of atomic and molecular systems; we refer the reader to
the survey \cite{Dehesa:01} and to references therein. Both
$S(\rho)$ and $ S(\gamma)$ also play role in a generalization of
the Heisenberg uncertainty relation: it has been established
\cite{BB:75} that for any pair of densities $\rho(\vec{r})$ and
$\gamma(\vec{p})$ in $D$-dimensional space, we have the sharp
inequality
\begin{equation} \label{xpeur}
S(\rho) + S(\gamma) \geq D(1+\ln \pi) \,,
\end{equation}
which expresses quantitatively the impossibility of simultaneous
localization of a pair of observables with no common eigenstates.

It is well known that the wave function of many important systems,
such as $D$-dimensional harmonic oscillator and hydrogen atom, are
expressible in terms of families of orthogonal polynomials. It is
not surprising that, as it has been shown in \cite{dvaypp,yvad94},
the computation of the entropies $S(\rho)$ and $ S(\gamma)$
usually can be reduced to integrals involving these polynomials.

Let $\nu$  be a positive unit Borel measure on $\Delta:=[-1,1]$
and let
\begin{equation*}\label{polynomials}
p_n(x)=\gamma _n \prod_{j=1}^{n} \left ( x-\zeta_j^{(n)}   \right
)\,, \quad \gamma _n>0,\quad n\in\N,
\end{equation*}
denote the corresponding sequence of \emph{orthonormal}
polynomials such that
\begin{equation*} \label{enint}
\int  p_n(x) p_m(x) \,   d \nu (x) = \delta_{mn}, \quad m, n \in\N\,.
\end{equation*}
We define the {\em information entropy} of the polynomials $p_n(x)$ as
\begin{equation} \label{spg}
E_n=E_n(\nu ) = - \int p^2_n(x)  \, \ln \big( p^2_n(x) \big) \,  d
\nu (x) \,.
\end{equation}
Throughout the paper, we will assume that the orthogonality measure
$\nu$ is absolutely continuous with respect to the Lebesgue
measure $\lambda$ on $\Delta$ with Radon--Nikodym derivative
$$d\nu/d\lambda=\nu'(x)=w(x),\quad w\in L^1(\Delta).$$
For
normalization purposes, we will always assume that the weight $w$
is unitary, i.e.
\begin{equation}\label{normal}
\int_\Delta w(x)dx=1\,.
\end{equation}
The information entropy will be
indistinctly denoted by $E_n(\nu)$ and $E_n(w)$. We follow this
convention below for other notations.

The asymptotic behavior of $E_n$ as $n \to \infty$ has a special
interest in the study of the so-called Rydberg states of
quantum-mechanical systems. Besides physical motivations, there
are some fascinating aspects of this problem because of a certain
universal behavior of related integrals, and because of a close
connection of the entropy $E_n$ with important functionals of the
normalized zero counting measures of the polynomials $p_n$,
\begin{equation*}
\label{zerocounting}
\mu_n =\frac{1}{n}\,  \sum_{j=1}^n \delta_{\zeta_j^{(n)}},\quad
n>0,
\end{equation*}
and of the following probability measures $\nu_n$ :
\begin{equation*}
\label{nu_n}
d \nu_n(x)=p^2_n(x)\, d \nu(x),\quad n\geq 0,
\end{equation*}
(note that $\nu_0=\nu$). Both measures are standard objects of
study in the analytic theory of orthogonal polynomials. For
instance, the normalized zero counting measure $\mu _n$ is closely
connected with the $n$-th root asymptotics of $p_n$, and as was
shown in \cite{Rakhmanov:77}, $\nu_n$ is associated with the
behavior of the ratio $p_{n+1}/p_n$ as $n \to \infty$.
\\[\baselineskip]
If $\mu$ and $\nu$ are positive Borel measures on $\C$, then their
{\em mutual entropy} is defined as
$$
S(\mu , \nu )=\begin{cases} -\infty, & \text{if $\mu $ is not $\nu
$-absolutely continuous,} \\
 -\int \log \left( \frac{d \mu }{d\nu } \right)\, d\mu , & \text{if
   $\mu $ is $\nu$-absolutely continuous,}
\end{cases}
$$
and their {\em mutual logarithmic energy} as
$$ I(\nu, \mu) = - \iint \ln |z-t|\,
d \nu (t) \, d\mu(z).
$$
With these notations
the entropy (\ref{spg}) is equivalently rewritten
as
\begin{equation}\label{entopy1}
E_n(\nu) = S( \nu_n   \, , \,   \nu  )=- 2 \ln \gamma _n + 2n \,
I(\mu _n, \nu _n)\,.
\end{equation}
In particular, from classical Jensen's inequality for integrals, it
follows immediately
that if both $\mu$ and $\nu$ are unit measures on $\Delta$, then
$S(\mu , \nu ) \leq 0$, with equality if and only if $\mu =\nu $. Hence,
$$
E_n(\nu) \leq 0 \, ,
$$
with equality if and only if $n=0$.

Aptekarev et al.\ \cite{Aptekarev:95} considered two subfamilies
of the usual Szeg\H{o} class of weights on $\Delta$, namely the Jacobi
weights and the Bernstein-Szeg\H{o} class (weights being bounded
above, bounded away from zero, and satisfying a Dini-Lipschitz
condition). In this last case it is known that the asymptotic
formula for the orthogonal polynomials $p_n$ holds uniformly in
$\Delta $, as $n$ tends to infinity. With these assumptions it has
been proved in \cite{Aptekarev:95}  that
\begin{equation} \label{sasha_limit}
   \lim_{n\to \infty} E_n(w) = S(\rho  \, , \, w  ) +
   \log(2) - 1 \,,
\end{equation}
where
$$\rho(x)=1/(\pi\sqrt{1-x^2})$$
denotes the Chebyshev unit weight on $\Delta$.
We are  concerned here with the problem whether a weaker form of
this equality holds in the Szeg\H o class of weights. We will
show that the right-hand side of (\ref{sasha_limit}) is actually
an asymptotic upper bound for the entropy $E_n(w)$ when
the weight $w$ satisfies the Szeg\H{o} condition (see assumption
(\ref{szego}) below). Furthermore, the expression (\ref{spg}) for
the entropy can be naturally split into two functionals, which
have simple asymptotic behaviors when $w$ belongs to the Bernstein
class (see Proposition \ref{prop:bernstein}). The situation with
the Bernstein class is in a sense optimal: the corresponding
limits provide asymptotic upper bounds for $w$ in the whole
Szeg\H{o} class. We also give conditions for the entropy and the
two functionals to tend to their upper bounds as the degree $n$
becomes large.

Finally, we must mention that in the case of an unbounded support
of the weight of orthogonality interesting results concerning the
asymptotics of the $E_n$ and related functionals have been
obtained recently in \cite{Levin03}.

\section{Statements of results}
The
weighted $L^p$ norm of a function $f$ with respect to a weight $k$ on
$\Delta$ will be denoted by
$$
\|f\|_{L^p(k)}=\left( \int_{\Delta} \left| f(x) \right|^p k(x)\,
dx \right)^{1/p}\,, \quad 1\leq p \leq\infty\,.
$$
We will simply write $L^p$ when $k \equiv 1$ on $\Delta$. \\
Though our main interest in this paper lies in the Szeg\H{o} class of
weights,
some other classes appear at different places. We recall
the definitions of these classes now.

The {\it Erd\H os-Turan class} $\EE\TT$ consists of
weights $w\in L^1$ such that $w>0$ almost everywhere on $\Delta$.

The {\it Szeg\H{o} class} $\SS$ consists of
weights $w\in L^1$ such that
\begin{equation}\label{szego}
    \log (w_0) \in{L^1(\rho)},
\end{equation}
where
$$w_0(x):=w(x)/\rho(x)=\pi \sqrt{1-x^2} w(x)$$
denotes the trigonometric weight corresponding to $w$.
The fact that $w_0\in L^1(\rho)$ implies
$\log^+ (w_0) \in {L^1(\rho)}$,
where, as usual, we denote
$$
\log^+(x)=\max \{ \log(x), 0\}\,, \quad x >0.
$$
Hence, condition (\ref{szego}) is
actually equivalent to
\begin{equation}\label{szego1}
S(\rho , w)=\int_\Delta\log(w_0(x))\rho(x)dx
>-\infty.
\end{equation}
Note that (\ref{szego}) and (\ref{szego1}) can equivalently be
rewritten as $\log (w)
\in{L^1(\rho)}$ and
$$\int_\Delta\log(w(x))\rho(x)dx
>-\infty$$
 respectively.

Finally the {\it Bernstein class} $\BB$ consists of
weights $w$ such that $w_0$ is given by the reciprocal of
a positive polynomial on $\Delta$.
As it is well-known, the class $\BB$ is an important class
useful for establishing asymptotic properties in the Szeg\H{o} theory of
orthogonal polynomials. Obviously, one has the following inclusions
$\BB\subset\SS\subset\EE\TT$.

We will also use the notations:
\begin{equation}\label{fnDef}
   f_n(x):=p_n(x)\sqrt{w_0(x)}\,,
\end{equation}
and for $M>0$,
\begin{equation}\label{defM}
\Delta_n(M) := \{ x \in \Delta : |f_n(x)| \geq M \}\,.
\end{equation}

One of the main results of the paper is:
\begin{theorem} \label{thm:boundE}
Assume that the weight $w$ belongs to the Szeg\H o class $\SS$.
Then, for all $M>\sqrt{2}$,
\begin{equation} \label{limit3}
    E_n(w)
      = S(\rho,w)+\log(2)-1   - \int_{\Delta_n(M)} p_n^2(x) \log^+(p^2_n(x) )\, w (x)\,
    dx  + o(1)\,, \quad n\to \infty\,.
\end{equation}
\end{theorem}
As a simple consequence of the above formula, we obtain the following
asymptotic upper bound together with necessary and sufficient conditions
for equality.

\begin{corollary}\label{cor:boundE}
Assume that the weight $w$ belongs to the Szeg\H o class $\SS$.
Then the following asymptotic upper bound for the entropy holds :
\begin{equation}\label{upperboundE}
\limsup_{n \to \infty} E_n(w)\leq S(\rho \, , \, w)+ \log(2) - 1
\,.
\end{equation}
Moreover, for a subsequence $n\in\Lambda\subset\N$,
\begin{equation}\label{limE}
\lim_{n \in \Lambda} E_n(w)= S(\rho \, , \, w)+\log(2) - 1 \,,
\end{equation}
if and only if there exists a constant $M>\sqrt{2}$, such that
\begin{equation}\label{suff_conditionE}
\lim_{n\in\Lambda}
    \int_{\Delta_n(M)} p_n^2(x) \log^+(p^2_n(x) )\, w (x)\,
    dx = 0.
\end{equation}
In this case (\ref{suff_conditionE}) is valid for for all
$M>\sqrt{2}$.

Furthermore, (\ref{suff_conditionE}) holds if there exists
$\varepsilon
>0$ such that either
\begin{equation}\label{suffCond2E}
  \sup_{n\in \Lambda} \int_\Delta \left( \log^+(p_n^2(x))
  \right)^{1+\varepsilon } p_n^2(x)w(x) \, dx < \infty \,
  \quad \mbox{or} \quad
  \sup_{n\in
  \Lambda} \int_\Delta \left( p_n^2(x) \right)^{1+\varepsilon } w(x)
  \, dx < \infty \,.
\end{equation}
\end{corollary}
\begin{remark}
Notice that the findings of \cite{Aptekarev:95} on
Bernstein-Szeg\H{o} polynomials are included in Corollary
\ref{cor:boundE} since for $w\in \BB$, $\log(w_0)$ is bounded and the
$f_n$ are uniformly bounded in $[-1,1]$. In contrast, the case of
Jacobi polynomials requires some extra considerations. One knows
that for the orthonormal Jacobi polynomials there exists a
constant $c$ such that for $n\geq 0 $ and $x \in [-1,1]$, $$
\left| P_n^{(\alpha,\beta)}(x) \right| \cdot \left(\sqrt{1-x} +
\frac{1}{n}
     \right)^{\alpha+1/2}\left(\sqrt{1+x} + \frac{1}{n}
     \right)^{\beta+1/2} \leq c/\sqrt{\pi}\, .
$$ Taking into account that here $w_0(x)=\pi
\left(1-x\right)^{\alpha+1/2}\left(1+x\right)^{\beta+1/2}$, we
find that for $p_n=P_n^{(\alpha,\beta)}$, $$
\left(p_n(x)\right)^{2+\varepsilon } w_0 (x) \leq c \sqrt{\pi} \,
\left(\frac{1-x}{\left(\sqrt{1-x} + 1/n \right)^{2+\varepsilon
}}\right)^{\alpha+1/2}\, \left(\frac{1+x}{\left(\sqrt{1+x} + 1/n
\right)^{2+\varepsilon }}\right)^{\beta+1/2}\,, $$ and the second
condition in (\ref{suffCond2E}) is satisfied.
\end{remark}
\begin{remark}
An inequality weaker than (\ref{upperboundE}) is a straightforward
consequence of the asymptotic behavior of the measures $\nu_n$.
Indeed, if $w\in\EE\TT$, we know from
Rakhmanov's Theorem \cite{Rakhmanov:77} that $d\nu_n(x) \to
\rho(x)dx$ as $n \to \infty$ in the weak-* topology. It follows
from the weak upper semicontinuity of the mutual entropy
\cite[Corollary~5.3]{Killip03} that $\limsup E_n(w)=\limsup
S(\nu_n,w) \leq S(\rho,w)$. In particular, it shows that if the
weight $w$ is in $\EE\TT\setminus\SS$,
$$
\lim_{n\to \infty} E_n(w)=-\infty\,.
$$
Nevertheless, it seems that a semicontinuity argument for the
entropy does not allow to explain the additional term $\log(2)-1$
occurring on the right-hand side of (\ref{limE}).
\end{remark}

The information entropy for Chebyshev polynomials orthonormal with
respect to $\rho$ has been computed in \cite{Dehesa:97,Yanez:99}:
\begin{equation}\label{caseChebyshev}
E_n(\rho)
= \log(2) - 1  \quad \text{for } n \geq 1 \,.
\end{equation}
Intuitively, Chebyshev polynomials are the most ``uniformly''
distributed polynomials, both for each $n$ and asymptotically as
$n \to \infty$. This fact is formally set in the next corollary.
\begin{corollary} \label{cor:chebyshev} If
\begin{equation}\label{larger}
\limsup_{n \to \infty} E_n(w)\geq  \log(2) - 1 \,,
\end{equation}
then $w=\rho$ and $E_n(w)=  \log(2) - 1$, $n\geq 1$.
\end{corollary}
The proof is a simple consequence of inequality
(\ref{upperboundE}). Indeed, from this inequality, we see that
(\ref{larger}) can only happen if $S(\rho,w)=0$ that is $\rho=w$.

\medskip

Now we exploit the connection between the entropy $E_n(w)$ and the
mutual energy $I(\mu_n,\nu_n)$ given in (\ref{entopy1}). It is
well known that in the class $\EE\TT$ both $\mu_n$ and $\nu_n$
tend (as $n \to \infty$) to the Chebyshev (equilibrium)
distribution given by the weight $\rho$ on $\Delta$. In
particular, from the convexity properties of the mutual energy it
follows that
$$
\lim_{n\to \infty} I(\mu_n,\nu_n)=I(\rho , \rho )=\log(2)\,.
$$
What is more surprising is that the next term of the asymptotic
expansion of $I(\mu_n,\nu_n)$ also exhibits a ``universal''
behavior, in the sense that it does not depend on the choice of
the weight $w$. Namely, if the entropy $E_n(w)$ satisfies
(\ref{limE}), then the following result is a direct consequence of
(\ref{entopy1}) and the well known asymptotic behavior of the
leading coefficient of $p_n$ (see (\ref{leading_coefficient})
below):
\begin{corollary}
\label{cor:energy} Assume $w$ is a weight in the Szeg\H o class
$\SS$ and condition (\ref{suff_conditionE}) is satisfied. Then the
mutual energy $I(\mu _n, \nu _n)$ has the following asymptotic
expansion:
$$
I(\mu _n, \nu
_n)=\log(2)-\frac{1}{2n}+o\left(\frac{1}{n}\right)\,,\quad
n\in\Lambda,\quad n \to \infty\,.
$$
\end{corollary}
This remarkable fact certainly deserves further study.

\medskip

Another aim of the paper is to study two related functionals $F_n$
and $G_n$, whose sum equals the entropy,
$$ E_n(w)=F_n(w)+G_n(w)\,, $$ and which are defined by
\begin{align}\label{Fn}
F_n(w) &= -\int_\Delta \log(p^2_n(x) w_0(x) )\, p^2_n(x) w(x) \,
dx= S(f_n^2 \rho \, , \, \rho ) \,, \\ \intertext{and}
\label{Gnbis} G_n(w) &= \int_\Delta \log( w_0(x) )\, p^2_n(x) w(x)
\, dx= -S(p^2_n w \, , \, p^2_n \rho)\,.
\end{align}

We will see that the functional $F_n$ also exhibits a
``universal'' behavior, while $G_n$ is sensitive to a particular
choice of the weight $w$, and is related naturally with the mutual
entropy $S(\rho,w)$.
Functionals $F_n$ and $G_n$ have a
particularly nice behavior for $w$ in the Bernstein class $\BB$:
\begin{proposition} \label{prop:bernstein}
Let $S$ be a polynomial of degree $2N$ ($N \geq 0$) such that
$S(x)>0$ for $x \in \Delta $, and assume that the orthogonality
weight satisfies
\begin{equation*}\label{bernstein}
    w_0(x)=\frac{1}{S(x)}\,, \quad x \in \Delta\, .
\end{equation*}
Then
\begin{equation}\label{caseBernstein}
F_n(w)= \log(2) - 1  \quad \text{for } n >N\,.
\end{equation}
Moreover,
\begin{equation}\label{limitG}
    \lim_{n \to \infty} G_n(w)=S(\rho \, , \, w)\,,
\end{equation}
and this limit takes place with a geometric rate. Consequently, the
same holds true for the limit in
(\ref{sasha_limit}).
\end{proposition}
The conjecture that constant entropy $E_n(w)$ is a (yet another)
characterization of Chebyshev polynomials (cf.\
(\ref{caseChebyshev})) belongs to L.\ Golinsky. We were able to
prove it in the Bernstein class $\mathcal B$.
\begin{proposition} \label{prop:cheb}
Let $w \in \mathcal B$ such that
$E_n(w)$ is constant for all sufficiently large $n$. Then $w=\rho
$.
\end{proposition}
Since Bernstein weights are suitable as approximation tool for the
whole Szeg\H{o} class, we could expect the asymptotic behavior
from Proposition \ref{prop:bernstein} to hold in a more general
setting. Nevertheless, the behavior of the entropy, as well as the
behavior of the two functionals $F_n$ and $G_n$, is extremely
sensitive to the growth of $p_n^2 w$, which may affect
convergence. In general, the following expression for the first
functional $F_n$ holds true:
\begin{theorem} \label{thm:boundF}
Assume the weight $w$ belongs to the Szeg\H o class $\SS$. Then, for all
$M>\sqrt{2}$,
\begin{equation}\label{identityForF}
 F_n(w)= \log(2) - 1   - \int_{\Delta_n(M)} \log(f_n^2(x)) \, f_n^2(x) \, \rho(x)
   \,dx+o(1)\,, \quad n \to \infty\,.
\end{equation}
\end{theorem}
Again, as a simple consequence of the above formula, we get the
\begin{corollary}
\label{cor:boundF}
Assume the weight $w$ belongs to the Szeg\H o class $\SS$. Then, the
following asymptotic upper bound for $F_n$ holds:
\begin{equation}\label{upperbound}
\limsup_{n \to \infty} F_n(w)\leq \log(2) - 1 \,.
\end{equation}
Moreover, for a subsequence $n\in\Lambda\subset\N$,
\begin{equation}\label{limF}
\lim_{n \in \Lambda} F_n(w)= \log(2) - 1 \,,
\end{equation}
if and only if there exists a constant $M>\sqrt{2}$, such that
\begin{equation}\label{suff_condition}
\lim_{n\in\Lambda} \int_{\Delta_n(M)} f_n^2(x) \log(f^2_n(x) )\,
\rho (x)\,     dx = 0\,,
\end{equation}
for $f_n$ and $\Delta_n(M)$ defined in (\ref{fnDef}) and
(\ref{defM}), respectively.
In this case, (\ref{suff_condition}) is
valid for every $M>\sqrt{2}$.

Furthermore, (\ref{suff_condition}) holds  if there exists an
$\varepsilon>0$ such that either
\begin{equation}\label{suffCond2}
  \sup_{n\in \Lambda} \int_\Delta \left( \log^+(f_n^2(x))
  \right)^{1+\varepsilon } f_n^2(x)  \, \rho(x)  dx < \infty \,
  \quad \mbox{or} \quad
  \sup_{n\in
  \Lambda} \int_\Delta \left(
  \, f_n^2(x) \right)^{1+\varepsilon } \, \rho(x)dx < \infty \,.
\end{equation}
\end{corollary}
\begin{remark} \label{remark:5}
The method of proof of Theorem \ref{thm:boundF} can be applied to
larger classes of weights. In fact, we only need an $L^2$
asymptotics of the polynomials $p_n$ on the support $\Delta$ of
the measure $\nu$, and that has been extended beyond the Szeg\H o
class. For instance, using our technique we can prove that
(\ref{upperboundE}) is valid for weights $w\in\mathcal{F}(dini)$,
introduced in \cite{Levin01}.
\end{remark}

\begin{remark}
Apparently, a necessary condition for (\ref{suffCond2}) is that
$w_0 \log (w_0) \in L^1(\rho)$ (cf.\ with  (\ref{szego})). If
$\log (w_0) \in L^\infty$ then there is equivalence between
conditions (\ref{suff_conditionE}) and (\ref{suff_condition}), and
between (\ref{suffCond2E}) and (\ref{suffCond2}), respectively.
\end{remark}
Concerning the second functional $G_n$, we use a result from
\cite{MNT} to deduce the following proposition.
\begin{proposition} \label{thm:boundG}
Assume the weight $w$
belongs to the Szeg\H{o} class $\SS$ and $\log^+(w_0)\in L^\infty$, then
\begin{equation}\label{secondLimSup}
\limsup_{n \to \infty} G_n(w)\leq S(\rho \, , \, w)=\int_\Delta
\log( w_0(x) )\, \rho(x) \, dx \,.
\end{equation}
Similarly, assume that $\log^-(w_0)\in L^\infty$, then
\begin{equation}\label{secondLimInf}
\liminf_{n \to \infty} G_n(w)\geq S(\rho \, , \, w)=\int_\Delta
\log( w_0(x) )\, \rho(x) \, dx \,.
\end{equation}
Hence, if $\log (w_0) \in L^\infty$, then
$$ \lim_{n \to
\infty} G_n(w)= S(\rho \, , \, w) \,.
$$
Furthermore, if the weight $w$ belongs to the set $\EE\TT\setminus\SS$, the
assumption
$\log^+(w_0)\in L^\infty$ still implies inequality
(\ref{secondLimSup}). In this case, (\ref{secondLimSup})
simplifies to $\lim_{n \to \infty} G_n(w)= -\infty$.
\end{proposition}
\section{Proofs of Theorems \ref{thm:boundE} and \ref{thm:boundF},
Corollaries \ref{cor:boundE} and \ref{cor:boundF}}
Before entering the proofs of our results, let us state two
preliminary lemmas. The first one is borrowed from \cite{Aptekarev:95}.
\begin{lemma}{\rm \cite[Lemma 2.1]{Aptekarev:95}}
\label{lemmaABD} Let $g$ be a continuous function on $\R$,
$g(\theta+\pi)=g(\theta)$, $f\in L^1([0,\pi])$, and let
$\gamma(\theta)$ be a function that is mesurable and almost
everywhere finite on $[0,\pi]$. Then, as $n\to\infty$,
$$\int_0^\pi g(n\theta+\gamma(\theta))f(\theta)d\theta\to
\frac{1}{\pi}\int_0^\pi g(\theta)d\theta \int_0^\pi
f(\theta)d\theta. $$
\end{lemma}
As remarked in \cite{Aptekarev:95}, when $\gamma(\theta)=0$ and
$g\in L^\infty[0,\pi]$, the statement of the lemma becomes a
well-known result of Fejer, cf. \cite[Chapter I, \S 20]{BAR}.

As second main ingredient in our proofs let us recall the
Szeg\H{o} asymptotics for $f_n(x) = \sqrt{w_0(x)}p_n(x)$: if
\begin{equation*}\label{defFandG}
g_n(x)=\sqrt{2}\, \cos(n \arccos x + \gamma(x))\,,
\end{equation*}
where
$$
\gamma (x)=\frac{1}{2\pi}\, \int_{\Delta} \frac{\log w_0(x)-\log
w_0(t)}{x-t} \, \sqrt{\frac{1-x^2}{1-t^2}}\, dt
$$
is the harmonic conjugate function to $\log w_0$, then in the
Szeg\H o class $\SS$, one has
\begin{equation}\label{limitSzego}
  \lim_{n \to \infty}  \| f_n-g_n\|_{L^2(\rho )}=0\,,
\end{equation}
and
\begin{equation} \label{leading_coefficient}
  \lim_{n \to \infty}
  \log \left( \frac{\gamma_n}{2^n} \right)=
  -\frac{1}{2}\, \left( \log(2)  +
  S(\rho  \, , \, w  )\right)\,.
\end{equation}
The mutual entropy in the right hand side of
(\ref{leading_coefficient}) is known as the Szeg\H{o} constant for
the weight $w$. Since the entropy integral is very sensitive to
the growth of $f_n^2=p_n^2 w_0$, the following lemma will be
useful; roughly speaking, it shows that the subsets $\Delta_n(M)$,
defined in (\ref{defM}), have no influence on the $L^2$
asymptotics (\ref{limitSzego}):
\begin{lemma}
For $w\in \SS$,
\begin{equation} \label{b1}
    \lim_{n \to \infty} \int_{\Delta_n(M)} \rho(x)\, dx = 0 .
\end{equation}
for every $M>\sqrt{2}$. Furthermore, let $ \widetilde f_n$, $n\geq
0$, be the sequence of truncated functions
\begin{equation}\label{defFtilde}
 \widetilde f_n(x) :=
    \left\{\begin{array}{ll}
        f_n(x) & \mbox{ for $x\in \Delta\setminus \Delta_n(M)$,} \\
        1 & \mbox{ for $x\in \Delta_n(M)$.} \end{array}\right.
\end{equation}
Then
\begin{equation} \label{b2}
    \lim_{n \to \infty}  \| \widetilde f_n-g_n\|_{L^2(\rho )}=0\,.
\end{equation}
\end{lemma}
\begin{proof}
 Observe first that by Cauchy-Schwarz inequality,
\begin{equation*}
   \int_{\Delta} | f_n^2(x)- g_n^2(x) | \rho(x)\, dx
    \leq
    \| f_n+g_n\|_{L^2(\rho )} \cdot \| f_n-g_n\|_{L^2(\rho )} \leq
    \left(\| f_n\|_{L^2(\rho )} +\|g_n\|_{L^2(\rho )}\right) \cdot \| f_n-g_n\|_{L^2(\rho
    )}\,,
\end{equation*}
so that
\begin{equation}\label{auxInequality}
\int_{\Delta} | f_n^2(x)- g_n^2(x) | \rho(x)\, dx
    \leq (1+\sqrt{2})\| f_n-g_n\|_{L^2(\rho )} \,.
\end{equation}
Now we can show that the Chebyshev (and hence, Lebesgue) measure
of $\Delta_n(M)$ is asymptotically vanishing: by
(\ref{auxInequality}),
\begin{eqnarray*}
   (M^2-2) \int_{\Delta_n(M)} \rho(x)\, dx &\leq&
   \int_{\Delta_n(M)} ( f_n^2(x)- 2 ) \rho(x)\, dx \leq
   \int_{\Delta_n(M)} | f_n^2(x)- g_n^2(x) | \rho(x)\, dx
   \\&\leq& (1+\sqrt{2})\| f_n-g_n\|_{L^2(\rho )} ,
\end{eqnarray*}
the right-hand side tending to zero as $n\to \infty$ by
(\ref{limitSzego}); this proves (\ref{b1}). Moreover, since
$|\widetilde f_n(x)|=1$ and $|g_n(x)|\leq \sqrt{2}$ for $x\in
\Delta_n(M)$, we have by (\ref{auxInequality})
\[
\begin{split}
 \| \widetilde f_n-g_n\|_{L^2(\rho )}^2 & =
\int_{\Delta\setminus \Delta_n(M)} | f_n(x)-g_n(x)|^2 \, \rho(x)
      \, dx+ \int_{ \Delta_n(M)} |\widetilde f_n(x)-g_n(x)|^2 \, \rho(x)
      \, dx \\ & \leq
      (1+\sqrt{2})\| f_n-g_n\|_{L^2(\rho )} + 3 \int_{\Delta_n(M)} \rho(x) \,
      dx\,.
\end{split}
\]
It remains to use (\ref{limitSzego}) and (\ref{b1}) to see that
(\ref{b2}) is satisfied. \qed
\end{proof}

\subsection{Proof of Theorem \ref{thm:boundE}}

Fix arbitrary $M>\sqrt{2}$ and let $\Delta_n(M)$ and $\widetilde
f$ be as defined in (\ref{defM}) and (\ref{defFtilde}),
respectively. We write the entropy as
\begin{equation} \label{b_aim}
   E_n(w)=S(f_n^2\rho,w) = S(g_n^2\rho,w)
   + [ S(\widetilde f_n^2\rho,w) - S(g_n^2\rho,w)]
   + [ S(f_n^2\rho,w) - S(\widetilde f_n^2\rho,w) ]
   .
\end{equation}
In three steps let us prove that the first term on the right has
as limit the first three terms in the right-hand side of (\ref{thm:boundE}),
the
second term tends to $0$, and the third term is asymptotically
negative and related to the integral in
(\ref{thm:boundE}).

Let $$\mathcal R(y)= y^2 \log(y^2), \quad y \in\R . $$ From Lemma
\ref{lemmaABD} we get
\begin{align} \nonumber &
    \lim_{n \to \infty} S(g_n^2\rho,w)
    \\ \nonumber &= -
    \lim_{n \to \infty} \int_0^\pi
    \mathcal R(g_n(\cos(\theta)))\frac{d\theta}{\pi}
    + \lim_{n \to \infty}
    \int_0^\pi \log(w_0(\cos(\theta))) g_n^2(\cos(\theta)) \frac{d\theta}{\pi}
    \\\nonumber  &= -
    \int_0^\pi
    \mathcal R(\sqrt{2}\cos(\theta))\frac{d\theta}{\pi}
    +
    \int_0^\pi 2\cos^2(\theta) \frac{d\theta}{\pi}
    \int_0^\pi \log(w_0(\cos(\theta))) \frac{d\theta}{\pi}
    \\ & = E_1(\rho) + S(\rho,w)=\log (2)-1 + S(\rho,w) . \label{Rightlimit}
\end{align}
Hence the first term on the right-hand side of (\ref{b_aim}) has
the required limit. The second term in (\ref{b_aim}) can be
written as
\begin{equation}\label{integral2}
    S\left(\widetilde f_n^2\rho,w\right) - S\left(g_n^2\rho,w\right) =
    \int_\Delta \left[
    \mathcal R\left(\frac{\widetilde f_n(x)}{\sqrt{w_0(x)}}\right)
    -     \mathcal R\left(\frac{g_n(x)}{\sqrt{w_0(x)}}\right) \right]
    w(x) \, dx .
\end{equation}
Recall that both $\widetilde f_n$ and $g_n$ are uniformly bounded
on $\Delta$ by $M$, and hence for $x\in \Delta$
$$
   \left|\mathcal R\left(\frac{\widetilde f_n(x)}{\sqrt{w_0(x)}}\right) \right|
    w_0(x) \leq |\mathcal R(\widetilde f_n(x))| + |\log(w_0(x))|\,
    \widetilde f_n^2(x)
    \leq M^2\log M^2 + M^2 |\log(w_0(x))| =: h(x) ,
$$ where $h\in L^1(\rho)$ by assumption (\ref{szego}). Similarly,
$$
     \left|\mathcal R\left(\frac{g_n(x)}{\sqrt{w_0(x)}}\right) \right| w_0(x) \leq h(x)\,,\qquad
     x\in \Delta\,.
$$
The integral in (\ref{integral2}) will be split into two parts
depending on whether $w_0$ is small or large. Fix an arbitrary $0<
\varepsilon< 1$; by the monotone convergence theorem there exists
a constant $C=C(\varepsilon )$ such that
$$
         0 \leq \int_{h(x)>C} h(x) \, \rho(x)\, dx =
         \int_{\Delta} h(x) \, \rho(x)\, dx -
         \int_{h(x)\leq C} h(x) \, \rho(x)\, dx
         < \varepsilon .
$$
Defining $\tau:=M^2\exp(-C/M^2)$ we see that $w_0(x)<\tau$ implies
that $h(x)>C$, and hence
$$
\left| \int_{w_0(x)<\tau} \left[
    \mathcal R\left(\frac{\widetilde f_n(x)}{\sqrt{w_0(x)}}\right)
    -     \mathcal R\left(\frac{g_n(x)}{\sqrt{w_0(x)}}\right) \right]
    w(x) dx \right|
    \leq 2 \int_{w_0(x)<\tau} h(x) \rho(x)\, dx
     \leq 2 \varepsilon.
$$
On the other hand, if $w_0(x)\geq \tau$, then
$$
\left|
    \frac{\widetilde f_n(x)}{\sqrt{w_0(x)}}\right|\leq
    \frac{M}{\sqrt{w_0(x)}}\leq  \frac{M}{ \sqrt{\tau} }=e^{C/(2 M^2) }=:C_1 ,
$$
and the same inequality is valid for $ g_n/\sqrt{w_0}$. Taking
into account that $\mathcal R $ is smooth,
\begin{eqnarray*} &&
    \left|
    \mathcal R\left(\frac{\widetilde f_n(x)}{\sqrt{w_0(x)}}\right)
    -     \mathcal R\left(\frac{g_n(x)}{\sqrt{w_0(x)}}\right)  \right|
    \leq \max_{|y| \leq C_1} | \mathcal R'(y) | \,
    \left|
    \frac{\widetilde f_n(x)}{\sqrt{w_0(x)}}
    -     \frac{g_n(x)}{\sqrt{w_0(x)}}  \right|
    \\&& \leq \max_{|y| \leq C_1}
          |  2 y (1+\log(y^2)) | \,
    \left|
    \frac{\widetilde f_n(x)}{\sqrt{w_0(x)}}
    -     \frac{g_n(x)}{\sqrt{w_0(x)}}  \right| \\
    && \leq C_2 \,
    \left|
    \frac{\widetilde f_n(x)}{\sqrt{w_0(x)}}
    -     \frac{g_n(x)}{\sqrt{w_0(x)}}  \right|\,,
\end{eqnarray*}
with $C_2:=\max \{4e^{-3/2}, 2 C_1 (1+\log(C_1^2)) \}$. Hence,
using the Cauchy-Schwarz inequality,
\[
\begin{split}
   \left| \int_{w_0(x)\geq\tau} \left[
    \mathcal R\left(\frac{\widetilde f_n(x)}{\sqrt{w_0(x)}}\right)
    -     \mathcal R\left(\frac{g_n(x)}{\sqrt{w_0(x)}}\right) \right]
    w(x) dx \right| &
    \leq C_2  \, \|  (\widetilde f_n - g_n)
    \sqrt{w_0(x)}\|_{L^1(\rho)}\\
   & \leq C_2  \, \|  \widetilde f_n - g_n
    \|_{L^2(\rho)},
\end{split}
\]
which by (\ref{b2}) tends to $0$ as $n \to \infty$. Taking into
account that $\varepsilon \in (0,1)$ was chosen arbitrary, we
conclude that
\begin{equation}\label{limit2}
    S\left(\widetilde f_n^2\rho,w\right) - S\left(g_n^2\rho,w\right)
    \longrightarrow 0\,, \quad n \to \infty\,.
\end{equation}

Thus, for establishing the expression for the entropy in
Theorem~\ref{thm:boundE},
it only remains to examine the last bracket on the right-hand side
of (\ref{b_aim}). Notice that since $\widetilde f_n=f_n$ on
$\Delta \setminus \Delta_n(M)$,
\begin{align*}
    S(f_n^2\rho,w) & - S(\widetilde f_n^2\rho,w)
    \\ = & - \int_{\Delta_n(M)} p_n^2(x) \log(p^2_n(x) )\, w (x)\,
    dx  + \int_{\Delta_n(M)}  \log\left(\frac{1}{w_0(x)}\right)\, \rho (x)\,
    dx
    \\ =&  - \int_{\Delta_n(M)} p_n^2(x) \log^+(p^2_n(x) )\, w (x)\,
    dx+ \int_{\widetilde \Delta_n(M)} p_n^2(x) |\log(p^2_n(x) )| \,w (x)\,
    dx
    \\ &    - \int_{\Delta_n(M)}  \log(w_0(x))\, \rho (x)\,
    dx\,,
\end{align*}
where
$$
\widetilde \Delta_n(M)=\left\{x\in\Delta_n(M):\, p_n^2(x)<1
\right\}\subset \Delta_n(M)\,.
$$
Observing that, for $p_n(x)^2\leq 1$, we have $0 \leq p_n^2(x)
|\log(p^2_n(x) )| \leq 1$, we obtain
$$
   0 \leq \int_{\widetilde \Delta_n(M)} p_n^2(x) |\log(p^2_n(x) )| \, w (x)\,
    dx \leq \int_{\Delta_n(M)} w (x)\, dx
    = \int_{\Delta_n(M)} w_0 (x)\,\rho(x)\, dx \,.
$$
Since $w_0\in L^1(\rho)$, $\log(w_0)\in L^1(\rho)$, by the
absolute continuity of the Lebesgue integral, relation (\ref{b1})
implies that
\begin{equation} \label{b3}
   \lim_{n\to \infty} \int_{\Delta_n(M)} w_0 (x)\,\rho(x)\, dx = 0 ,
   \quad \mbox{and} \quad
   \lim_{n\to \infty} \int_{\Delta_n(M)} \log(w_0 (x))\,\rho(x)\, dx = 0
   \,,
\end{equation}
showing that
\begin{equation}\label{auxIneq}
     S(f_n^2\rho,w)  - S(\widetilde f_n^2\rho,w)=- \int_{\Delta_n(M)} p_n^2(x) \log^+(p^2_n(x) )\, w (x)\,
    dx  + o(1)\,, \quad n\to \infty\,.
\end{equation}
Hence, gathering (\ref{Rightlimit}), (\ref{limit2}), and
(\ref{auxIneq}) in (\ref{b_aim}), we get (\ref{limit3}). \qed

\subsection{Proof of Corollary \ref{cor:boundE}}
Since $$
\int_{\Delta_n(M)} p_n^2(x) \log^+(p^2_n(x) )\, w (x)\,
    dx \geq 0\,,
 $$
relation (\ref{upperboundE}) is a trivial consequence of
Theorem~\ref{thm:boundE}. Suppose now that (\ref{suff_conditionE})
holds for some $M>\sqrt{2}$, then (\ref{limE}) follows immediately
from (\ref{limit3}). Conversely, if (\ref{limE}) is true then it
follows from Theorem~\ref{thm:boundE} that (\ref{suff_conditionE})
holds for all $M>\sqrt{2}$.

In order to prove that (\ref{suffCond2E}) is
sufficient for (\ref{suff_conditionE}), notice that, by
H\"older's inequality,
\begin{equation} \label{ineqSuffCond}
\begin{split}
     & \int_{\Delta_n(M)}   p_n^2(x) \log^+(p^2_n(x) )\, w (x)\,
    dx \\ & \leq
     \left(\int_{\Delta_n(M)} p_n^2(x) (\log^+(p^2_n(x) ))^{1+\epsilon}\, w (x)\, dx
     \right)^{\frac{1}{1+\epsilon}} \,
     \left(\int_{\Delta_n(M)} p_n^2(x) w (x)\, dx \right)^{1 - \frac{1}{1+\epsilon}}\,
     .
\end{split}
\end{equation}
Furthermore,
\begin{eqnarray*} &&
  \int_{\Delta_n(M)} p_n^2(x) w (x)\, dx
  \leq \int_{\Delta_n(M)} [ f_n^2(x)-g_n^2(x)] \rho (x)\, dx +
  \int_{\Delta_n(M)} g_n^2(x) \rho (x)\, dx \\ && \leq \int_{\Delta_n(M)} [ f_n^2(x)-g_n^2(x)] \rho (x)\, dx + 2
  \int_{\Delta_n(M)} \rho (x)\, dx
  \\&&
  \leq (1+\sqrt{2})\| f_n - g_n \|_{L^2(\rho)}
  + 2 \int_{\Delta_n} \rho (x)\, dx =o(1)\,, \quad n \to \infty\,,
\end{eqnarray*}
where we have used (\ref{limitSzego}), (\ref{b1}) and
(\ref{auxInequality}).

If we assume that the first condition in (\ref{suffCond2E}) holds,
then the first factor on the right hand side of
(\ref{ineqSuffCond}) is uniformly bounded  in $n$, and
(\ref{suff_conditionE}) follows.

Finally, notice that the second condition in (\ref{suffCond2E})
implies the first one since $\log^+(z)\leq z$ for $z\geq 0$, and
hence
$$
    (\log^+(y))^{1+\varepsilon} =
    \left(\frac{1+\varepsilon}{\varepsilon}\right)^{1+\varepsilon} \,
    (\log^+(y^{\frac{\varepsilon}{1+\varepsilon}}))^{1+\varepsilon}
    \leq   \left(\frac{1+\varepsilon}{\varepsilon}\right)^{1+\varepsilon} \,
    y^\varepsilon , \quad y \geq 0 .
$$
\qed

\subsection{Proof of Theorem \ref{thm:boundF}}

Our proof for Theorem~\ref{thm:boundF} follows closely the
arguments of the proof of Theorem~\ref{thm:boundE}, but some parts
simplify. As before let $\mathcal R(y)= y^2 \log(y^2)$, $y \in\R$,
and fix $M>\sqrt{2}$. We write the functional as follows:
\begin{eqnarray}\label{aim2}
   F_n(w)&=& \int_\Delta [-\mathcal R(g_n(x))] \, \rho(x)\, dx
   \\&&\nonumber
    + \int_\Delta [\mathcal R(g_n(x))-\mathcal R(\widetilde f_n(x))] \, \rho(x)\, dx
   + \int_\Delta [\mathcal R(\widetilde f_n(x))-\mathcal R(f_n(x))] \,
   \rho(x)\, dx.
\end{eqnarray}
Here the first integral on the right-hand side of (\ref{aim2}) has
the limit $E_1(\rho)=\log(2)-1$ by Lemma~\ref{lemmaABD}. The last
one can be written as $$
   \int_\Delta [\mathcal R(\widetilde f_n(x))-\mathcal R(f_n(x))] \,
   \rho(x)dx = - \int_{\Delta_n(M)} \log(f_n^2(x)) \, f_n^2(x) \, \rho(x)
   \,dx \leq 0 ,
$$ the right-hand side coinciding with the integral in (\ref{identityForF}).
Thus Theorem \ref{thm:boundF}
follows by showing that the second integral on the right-hand side
of (\ref{aim2}) is asymptotically vanishing. Recalling that
$|\widetilde f_n(x)|$ and $|g_n(x)|$ are uniformly bounded by $M$
for all $n\geq 0$ and $x\in \Delta$, we obtain
\begin{eqnarray*} &&
   \left| \int_\Delta [\mathcal R(g_n(x))-\mathcal R(\widetilde
   f_n(x))] \, \rho(x)dx
   \right|
   \leq \max_{y\in[-M,M]} | \mathcal R'(y) | \,
   \int_\Delta \left| g_n(x) - \widetilde f_n(x)) \right| \,
   \rho(x)\, dx
   \\&& \leq M^2(1+\log M^2) \, \| g_n - \widetilde f_n \|_{L^1(\rho)}
   \leq M^2(1+\log M^2) \, \| g_n - \widetilde f_n \|_{L^2(\rho)} ,
\end{eqnarray*}
the term on the right tending to zero as $n \to \infty$ by
(\ref{b2}). \qed

\subsection{Proof of Corollary \ref{cor:boundF}}
Since $$
\int_{\Delta_n(M)} f_n^2(x) \log(f^2_n(x) )\, \rho (x)\,
    dx \geq 0\,,
 $$
relation (\ref{upperbound}) is a trivial consequence of
Theorem \ref{thm:boundF}.
Suppose now that (\ref{suff_condition}) holds for some $M>\sqrt{2}$,
then (\ref{limF}) follows immediately from (\ref{identityForF}). Conversely,
if (\ref{limF}) is true then it follows from
Theorem~\ref{thm:boundF} that (\ref{suff_condition}) holds for
all $M>\sqrt{2}$.

In order to prove that the first condition in (\ref{suffCond2})
(which clearly is weaker than the second one)
is
sufficient for (\ref{suff_condition}), notice that, by
H\"older's inequality,
\begin{eqnarray*} &&
     \int_{\Delta_n} f_n^2(x) \log^+(f^2_n(x) )\, \rho (x)\,
    dx \\&& \leq
     \left(\int_{\Delta_n} f_n^2(x) (\log^+(f^2_n(x) ))^{1+\varepsilon}\, \rho (x)\, dx
     \right)^{\frac{1}{1+\varepsilon}} \,
     \left(\int_{\Delta_n} f_n^2(x) \rho (x)\, dx \right)^{1 - \frac{1}{1+\varepsilon}}\, ,
\end{eqnarray*}
and we may conclude as in the proof of Corollary \ref{cor:boundE}
that the second factor in the right hand side tends to zero. \qed

\section{Proofs of Propositions \ref{prop:bernstein}, \ref{prop:cheb},
  and
\ref{thm:boundG}}
\subsection{Proof of Proposition \ref{prop:bernstein}}

Let us make the change of variables $ x= (z+1/z )/2 $. It is well
known that since $S(x)>0$ on $\Delta$ we may write $S$ as
\begin{equation}
\label{def:q}
S(x)
=\left|q(z)\right|^2=q\left( z \right) q\left( 1/z \right)
\end{equation}
with
$q$ a polynomial of degree $2N$ with real coefficients having all
its zeros outside the disk and $q(0)>0$. Moreover
 \begin{equation} \label{bernsteinOP}
     p_n(x) = \frac{1}{\sqrt{2}}\, \left( z^n q\left(z^{-1}\right) + z^{-n}\, q(z) \right)
   \end{equation}
is the orthonormal polynomial of degree $n>N$ with respect to the
Bernstein weight $\rho/S$.
Introducing the Blaschke product
\begin{equation}
\label{Blaschke-def}
B_n(z)=z^{2n} q(1/z)/q(z), \quad n\geq N,
\end{equation}
we find that
   \begin{equation*} \label{bernstein_blaschke}
    p_n^2(x) w_0(x) =\frac{1}{2}\,
     | 1 + B_{n}(z)|^2
      =1 + \frac{1}{2} (B_{n}(z) +B_{n}(1 / z)) , \qquad |z|=1\,.
   \end{equation*}
Since, for $n>N$, $B_{n}(0)=0$, and $B_{n}$ is analytic in the disk, we have
   \begin{align*}
    \log(2) -F_n(w)&=  \log(2) + \int \log\left(p_n^2(x) w_0(x) \right) p_n^2(x) w_0(x)
     \rho(x)\, dx
     \\&
     =\frac{1}{2\pi} \int_{|z|=1} \log\left(| 1 + B_{n}(z)|^2\right) \left[ 1 + \frac{1}{2} (B_{n}(z) +
     B_{n}(1 / z)) \right]
     \, |dz|
     \\&
     = \Re \left(
     \frac{1}{2\pi i} \int_{|z|=1} \log( 1 + B_{n}(z)) \, [ 2 +
     B_{n}(z) + B_{n}(1 / z) ]\,  \frac{dz}{z} \right)\,.
   \end{align*}
Since $|B_n(z)|<1$ for $|z|<1$, the function $\log( 1 + B_{n}) \, [ 2 +
B_{n}]$ is holomorphic inside the disk and vanishes at the origin.
Thus,
\begin{align}\label{comput-F}
  \log(2) -F_n(w)&= \Re \left(
     \frac{1}{2\pi i} \int_{|z|=1} \log( 1 + B_{n}(z)) \, B_{n}(1 / z) \,  \frac{dz}{z} \right)\\
     &= \Re \left(\frac{1}{2\pi i}\, \int_{|z|=1} \frac{\log( 1 +
    B_{n}(z))}{B_{n}(z)}\, \frac{dz}{z}\right) \,,
\end{align}
where we have used that $B_n(1/z)=1/B_n(z)$. Observe that the last
integrand is analytic in a neighborhood of the unit circle, and we
can integrate along a smaller circle $|z|=r<1$, where
$|B_n(z)|<1$. Replacing $\log$ by its  uniformly convergent Taylor
expansion we get finally that this integral equals 1, which proves
(\ref{caseBernstein}).

On the other hand, by a similar reasoning we have
   \begin{align}
    G_n(w)&=   \int_{-1}^1 \log(w_0(x)) p_n(x)^2 w(x)dx
\nonumber \\&       = - 2 \Re\left(  \frac{1}{2\pi i} \int_{|z|=1}
\log(q(z)) [ 1 + \frac{1}{2} (B_{n}(z) + B_{n}(1 / z)) ]
      \, |dz|\right)
 \nonumber \\&
       = - 2 \log(q(0)) - \frac{1}{2\pi i}\, \int_{|z|=1} \log(q(z))
       B_{n}(1 / z)\,
       \frac{dz}{z}. \label{Blaschke}
   \end{align}
Note that in the last expression of (\ref{Blaschke}), taking the real part is
not necessary since $q$ and $B_n$ are real functions.
Integrating now along  $|z|=R>1$, we observe that
   $|B_n(1/z)|$ becomes there geometrically small, which yields a geometric
rate of convergence for
$$
\lim_{n\to \infty} G_n(w)= - 2 \log(q(0)) = -\Re \frac{1}{\pi}\,
\int_{|z|=1}\log \left( q(z)\right) \, \frac{dz}{z}= S(\rho \, , w)\,,
$$
which proves (\ref{limitG}). \qed

\subsection{Proof of Proposition \ref{prop:cheb}}
From the computations of $F_n(w)$ and $G_n(w)$ in the proof of Proposition
\ref{prop:bernstein}, see (\ref{comput-F}) and (\ref{Blaschke}), we know that
$E_n(w)$ is constant for $n$ large, say $n>N_0>N$, if and only if
\begin{equation}
\label{int-nul}
\frac{1}{2\pi i}\, \int_{|z|=1}
       \frac{B_N(1/z)\log(q(z))}
       {z^{2n-2N}}\,
       \frac{dz}{z}= 0, \quad n>N_0,
\end{equation}
where the polynomial $q$ and the Blaschke product $B_N$ are defined by
(\ref{def:q}) and (\ref{Blaschke-def}) respectively.
Since $\log(q(z))$ is analytic in some neighborhood $\UU$
of the unit disk, we may conclude that $\log(q(z))B_N(1/z)$ is
meromorphic in $\UU$, and thus can be
written as
\begin{equation}
\label{decomp}
      B_N(1/z)\log(q(z)) = r(z)  + f(z) , \qquad z\in\UU,
\end{equation}
where $r$ is a rational function such that $z^{2N}q(1/z)r(z)$ is a polynomial
of degree at most $2N-1$, and $f$ is
analytic in $\UU$.
Since $r$ is analytic outside the unit disk and grows like at most $1/z$ at
infinity, we deduce
$$
      \frac{1}{2\pi i}\, \int_{|z|=1}
       \frac{r(z)}
       {z^{2n-2N}}\,
       \frac{dz}{z}=0,\quad n>N_0,
$$
which implies, together with (\ref{int-nul}) and (\ref{decomp}), that
$$
             \frac{1}{2\pi i}\, \int_{|z|=1}
       \frac{f(z)}
       {z^{2n-2N}}\,
       \frac{dz}{z}=0\,,\quad n>N_0.
$$
Hence, all sufficiently high even Taylor coefficients of $f$ vanish. As a
consequence, $f(z)+f(-z)=P(z)$ is a
polynomial, and
\begin{equation}\label{analytic}
     B_N(1/z)\log(q(z))+ B_N(-1/z)\log(q(-z)) = r(z)+r(-z) + P(z), \quad
     |z|\leq 1 .
\end{equation}
Since the right-hand side of (\ref{analytic}) is a rational
function, the principle of analytic continuation applies, showing
that (\ref{analytic}) actually holds everywhere in $\C$. First,
assume that the polynomial $q$ is even, that is $q(z)=q(-z)$,
$z\in\C$. Then, it follows from (\ref{analytic}) that $\log(q(z))$
is a rational function so that $q$ can only be a constant, namely
1 by the normalization (\ref{normal}) of the weight $w$. Second,
assume that the polynomial $q$ is not even (hence different from a
constant). It implies the existence of some root $\alpha\in\C$ of
$q$ such that either $q(-\alpha)\neq 0$ or $-\alpha$ is a root of
$q$ of different multiplicity than that of $\alpha$. Note that
$\alpha\neq 0$ since, by assumption, $q(0)>0$. Then we get a
contradiction. Indeed, in view of the definition
(\ref{Blaschke-def}) of $B_N$, we readily observe that the
left-hand side of (\ref{analytic}) has a branch point at $\alpha$
while the right-hand side has not. Hence, $q(z)$ is constant,
equal to 1, and the proof of Proposition~\ref{prop:cheb} is
finished. \qed

\subsection{Proof of Proposition \ref{thm:boundG}}
Choosing $p=2$ and $g=|\log(w_0)w|^{1/2}\in L^1$ in Theorem 2 of
\cite{MNT} shows that
\begin{equation}\label{mnt}
\liminf_{n \to \infty} \int_\Delta|\log( w_0(x) )|p_n^2(x)w(x)dx
\geq \int_\Delta|\log( w_0(x) )| \rho(x) \, dx
\end{equation}
for any weight $w$ in the Erd\H os-Turan class $\EE\TT$. If $\log^+( w_0)\in
L^\infty$, there exists a constant $C>1$ such that $w_0(x)\leq C$,
$x\in\Delta$. Hence $|\log( w_0/C)|=-\log( w_0/C)$ and substracting
$\log(C)$ to both sides of
(\ref{mnt}), we get
(\ref{secondLimSup}) since
$$\int_\Delta{\rho(x)}dx=
\int_\Delta{p_n^2(x)w(x)}dx={1}.
$$
A similar reasonning shows (\ref{secondLimInf}) when $\log^-(
w_0)\in L^\infty$. Since this argument applies for any weight in
the Erd\H os-Turan class, the last assertion in the proposition
also follows. \qed

\section*{Acknowledgement} This work was partially supported by
INTAS project 2000--272 (B.B., A.M.F.\ and F.W.), a research grant
from the Ministry of Science and Technology (MCYT) of Spain,
project code BFM2001-3878-C02 (B.B and A.M.F.), by Junta de
Andaluc{\'\i}a, Grupo de Investigaci{\'o}n FQM 0229 and by the Ministry of
Education, Culture and Sports of Spain through the grant
PR2003--0104 (A.M.F.). A.M.F.~wishes to acknowledge also the
hospitality of the Department of Mathematics, Universite  de Lille
I,  France, where this work was started.

E.A.R.\ acknowledges also the support of a research grant from the
Ministry of Science and Technology (MCYT) of Spain, project code
SAB2001-0120, and the hospitality of the University of
Almer{\'\i}a.


\begin{thebibliography}{10}

\bibitem{Aptekarev:95}
A.~I. Aptekarev, V.~S. Buyarov, and J.~S. Dehesa.
\newblock Asymptotic behavior of the {$L^p$}-norms and the entropy for general
  orthogonal polynomials.
\newblock {\em Russian Acad.\ Sci.\ Sb.\ Math.}, 82(2):373--395, 1995.

\bibitem{BAR} N.\ K. Bari.
\newblock {\em A treatise on trigonometric series}, Pergamon Press,
Oxford, 1964.

\bibitem{BB:75}
I.~Bialynicki-Birula and J.~Mycielsky.
\newblock  Uncertainty Relations for Information Entropy in Wave
Mechanics.
\newblock {\em Commun.\ Math.\ Phys.}, 44:129--132, 1975.

\bibitem{Dehesa:97}
J.~S. Dehesa, W.~Van~Assche, and R.~J. Y{\'a}{\~n}ez.
\newblock Information entropy of classical orthogonal polynomials and their
  application to the harmonic oscillator and {C}oulomb potentials.
\newblock {\em Methods and Appl.\ Analysis}, 4:91--110, 1997.

\bibitem{Dehesa:01}
J.~S. Dehesa, A.~Mart{\'{\i}}nez-Finkelshtein, and
J.~S{\'a}nchez-Ruiz.
\newblock  Quantum Information Entropies and Orthogonal Polynomials.
\newblock {\em J. Comput. Appl. Math.}, 133:23--46, 2001.

\bibitem{dvaypp} J.S. Dehesa, W. Van Assche, R.J. Y\'{a}\~{n}ez,
Information entropy of classical orthogonal polynomials and their
application to the harmonic oscillator and Coulomb potentials,
Meth. Appl. Anal. 4 (1997) 91--110.

\bibitem{dg90} R.M. Dreizler and E.K.U. Gross.
\newblock \emph{Density Functional Theory: An Approach to the Quantum
Mechanics}.
 \newblock Springer-Verlag, Heidelberg, 1990.

\bibitem{Killip03}
R.~Killip and B.~Simon.
\newblock Sum rules for {J}acobi matrices and their applications to spectral
  theory.
\newblock {\em Annals of Math.},
158:253--321, 2003.

\bibitem{Levin01}
A.~L. Levin and D.~S. Lubinsky.
\newblock {\em Orthogonal Polynomials with Exponential Weights}, volume~4 of
  {\em CMS Books in Mathematics}.
\newblock Springer Verlag, 2001.

\bibitem{Levin03}
E.~Levin and D.~S. Lubinsky.
\newblock Asymptotics for entropy integrals associated with exponential
  weights.
\newblock {\em J.\ Comput.\ Appl.\ Math.}, 156:265--283, 2003.

\bibitem{MNT} A. Mat\'e, P. Nevai and V. Totik.
\newblock Necessary conditions for weighted mean convergence of
Fourier series in orthogonal polynomials.
\newblock {\em J. Approx. Theory}, 46:314--322, 1986.

\bibitem{py89} R.G. Parr, W. Yang,
 \newblock {\em Density Functional Theory of
Atoms and Molecules},
 \newblock Oxford University Press, New York,
1989.

\bibitem{Rakhmanov:77}
E.~A. Rakhmanov.
\newblock On the asymptotics of the ratio of orthogonal polynomials.
\newblock {\em Math. USSR Sb.}, 32:199--213, 1977.

\bibitem{sha48}  C.E. Shannon, A mathematical theory of communication,
Bell Syst. Tech. J. 27 (1948) 379-423, 623-656; reprinted in: The
Mathematical Theory of Communication, Eds. C. E. Shannon and W.
Weaver, University of Illinois Press, Urbana, 1949.

\bibitem{yvad94} R.J. Y\'{a}\~{n}ez, W. Van Assche, J.S. Dehesa, Position
and momentum information entropies of the D-dimensional harmonic
oscillator and hydrogen atom, Phys. Rev. A 50 (1994) 3065--3079.

\bibitem{Yanez:99}
R.~J. Ya{\~n}ez, W.~Van~Assche, R.~Gonz{\'a}lez-F{\'e}rez, and
J.~S. Dehesa.
\newblock Entropic integrals of hyperspherical harmonics and spatial entropy of
  {$D$}-dimensional central potentials.
\newblock {\em J. Math. Physics}, 40(11):5675--5686, 1999.

\end{thebibliography}
\end{document}